\documentclass[a4paper,11pt]{amsart}
\usepackage{extarrows}
\usepackage{amssymb}
\usepackage{mathrsfs}
\usepackage{amsmath,amssymb,amsthm,latexsym,amscd,mathrsfs}
\usepackage{indentfirst}
\usepackage{extarrows}
\usepackage[all]{xy}

\usepackage{multirow}

\newtheorem{theorem}{Theorem}
\newtheorem{problem}[theorem]{Problem}

\newtheorem{conjecture}{Conjecture}

\renewcommand{\leq}{\leqslant}
\renewcommand{\geq}{\geqslant}

 \makeatletter

\newcommand{\Diff}{\mathrm{Diff}}
\newcommand{\const}{\mathrm{const}}
\newcommand{\id}{\mathrm{id}}
\newcommand{\II}{\mathrm{II}}
\newcommand{\Ric}{\mathrm{Ric}}
\newcommand{\TangConj}{\mathbf{Tang~ Conjecture}}
\newcommand{\GroveConj}{\mathbf{Grove~ Conjecture}}
\newcommand{\ChernConj}{\mathbf{Chern~ Conjecture}}

\makeatother
\title[Problems related to isoparametric theory]{Problems related to isoparametric theory}
\author[J. Q. Ge]{Jianquan Ge}
\address{School of Mathematical Sciences, Beijing Normal
University, Beijing 100875, P. R. China}
\email{jqge@bnu.edu.cn}

\date{}
\thanks{The author is supported by Beijing Natural Science Foundation (Z190003), NSFC (No. 11522103, 11331002) and by Fundamental Research Funds for the Central Universities of China.}
\begin{document}
\maketitle


\begin{abstract}
 In this note we briefly survey and propose some open problems related to isoparametric theory.
\end{abstract}

\section*{1. Introduction}

In recent years the theory of isoparametric hypersurfaces has been extensively studied with fruitful results and interesting applications.  Although we have arrived at a celebrated ending of the long way of classification for isoparametric hypersurfaces in unit spheres since from E. Cartan in the late 1930s', we find that we are standing at a new beginning of exploring applications of this theory as well as an even longer way of classification for isoparametric hypersurfaces in general Riemannian manifolds such as homotopy spheres. For details we refer to the excellent book \cite{CR15} and the surveys \cite{Chi17, QG16, QT18, TY17}.

For the sake of convenience, we briefly recall some notations. A hypersurface of a Riemannian manifold is called \textit{isoparametric} if its nearby parallel hypersurfaces have constant mean curvature; or equivalently, it is locally a regular level set of an \textit{isoparametric function} $f$, by definition,
\begin{equation}\label{isop-def}
\left\{\begin{array}{ll}
|\nabla f|^2=b(f), &\\
\Delta f=a(f),
  \end{array}
  \right.
\end{equation}
where $b,a$ are $C^2$ and $C^1$ functions on $\mathbb{R}$; or a regular leaf of an \textit{isoparametric foliation} (i.e., a singular Riemannian foliation of codimension $1$ with constant mean curvature regular leaves).
A function with only the first equation of (\ref{isop-def}) satisfied is called a \textit{transnormal function} whose level sets constitute a singular Riemannnian foliation of codimension $1$.  In the unit sphere $\mathbb{S}^{n+1}$ (or real space forms), E. Cartan showed that a hypersurface is isoparametric if and only if it has constant principal curvatures. The number $g$ of distinct principal curvatures must be one of $\{1,2,3,4,6\}$ shown by M\"{u}nzner \cite{Mu80} with rather complicated topological proof (see Fang \cite{Fa17} for another proof using rational homotopy theory). The maximum and minimum level sets of $f$ are smooth submanifolds called the focal submanifolds, denoted by $M_{\pm}$ with codimensions $m_{\pm}+1$. In particular, isoparametric hypersurfaces $M^n$ in the unit sphere $\mathbb{S}^{n+1}$ correspond to the level sets (intersection with the sphere) of the so-called Cartan-M\"{u}nzner (isoparametric) polynomials $F(x)$ ($x\in\mathbb{R}^{n+2}$) which are homogeneous polynomials of degree $g$ satisfying
\begin{equation}\label{CMeq}
\left\{\begin{array}{ll}
|\nabla F|^2 =g^2|x|^{2g-2}, &\\
\Delta F=\frac{g^2}{2}(m_{-}-m_{+})|x|^{g-2},
  \end{array}
  \right. \quad x\in\mathbb{R}^{n+2}.
\end{equation}

Classically, classifying isoparametric hypersurfaces is done in the category of submanifold geometry and thus is carried out up to isometries of the ambient space. Two isoparametric hypersurfaces (foliations) are regarded as the same once they can be moved to one another by an isometry of the ambient space. This produces exactly the fundamental theory of isoparametric hypersurfaces in Riemannian manifolds with good symmetries such as the real space forms or Riemannian symmetric spaces. Moreover, in unit spheres $\mathbb{S}^{n+1}$ (as $n$ varies), there are infinitely many homogeneous and inhomogeneous isoparametric foliations, each with a unique minimal isoparametric hypersurface (with constant principal curvatures) where the simplest cases $g=1$ and $g=2$ are the totally geodesic sphere $\mathbb{S}^n$ and the minimal Clifford tori $\mathbb{S}^k(\sqrt{\frac{k}{n}})\times \mathbb{S}^{n-k}(\sqrt{\frac{n-k}{n}})$ ($1\leq k\leq [\frac{n}{2}]$). Therefore, plenty of applications of this fundamental theory have been explored (e.g., \cite{GX10, PT96, QT16, Sh00, PP06, MO14, Mi13, Sa16} and references in the book and surveys above), and even more are expected as we will see later in Section 2 for the classical isoparametric theory.

Based on Qian-Tang's fundamental construction (\cite{QT15}), classifying isoparametric foliations can be alternatively done in the category of differential topology which is carried out up to foliated diffeomorphisms of the ambient space. Two isoparametric foliations $(N_1, \mathfrak{F}_1)$ and $(N_2, \mathfrak{F}_2)$ are called equivalent if there is a foliated diffeomorphism between them, i.e., there is a diffeomorphism from $N_1$ to $N_2$ which maps each leaf of $\mathfrak{F}_1$ to that of $\mathfrak{F}_2$. For details we refer to the survey \cite{QG16}. Applications to other areas such as $4$ dimensional smooth Poincar\'{e} conjecture, $4$-manifolds geometry, closed geodesic and exotic smooth structures, etc., have been achieved by \cite{GT13, GR15, TZ14, Ge16}. More explorations are expected as we will see later in Section 3 for the general isoparametric theory.

\section*{2. Problems related to isoparametric theory on unit spheres}
\begin{problem}
$(\ChernConj)$ Let $M^n$ be a closed minimal hypersurface with constant scalar curvature (CSC) in the unit sphere $\mathbb{S}^{n+1}$. Then $M$ is isoparametric.
\end{problem}

This is a refined version of Chern's original problem (cf. \cite{CDK70, Yau82}): whether the set of scalar curvatures of such minimal CSC hypersurfaces in $\mathbb{S}^{n+1}$ is discrete. For isoparametric case, the set is included in $\{(n-g)n \mid g=1,2,3,4,6\}$.

This problem originated from Simons's pinching rigidity theorem and generated a large amount of studies on such pinching or gap rigidity (cf. \cite{GT12, XX17} and references therein). The Chern Conjecture is trivially true when $n\leq2$, and has been proven affirmatively when $n=3$ by Chang \cite{Chang93}, in which case the conjecture even holds for constant mean curvature CSC hypersurfaces (cf. \cite{DB90} and see \cite{TWY18} for a generalization). It is still widely open for the case $n\geq4$ in general.

\begin{problem}
Compute the spectrum of minimal isoparametric hypersurfaces in $\mathbb{S}^{n+1}$.
\end{problem}

The spectrum for the $g=1$ case ($\mathbb{S}^{n}$) and the $g=2$ case ($\mathbb{S}^k(\sqrt{\frac{k}{n}})\times \mathbb{S}^{n-k}(\sqrt{\frac{n-k}{n}})$) are well known. The $g=3$ case (tubes around the Veronese projective planes $\mathbb{FP}^2$ for $\mathbb{F}=\mathbb{R}, \mathbb{C}, \mathbb{H}, \mathbb{O}$) has been completely computed by Solomon \cite{So90}. For the infinite family in the case $g=4$ and the only two homogeneous minimal isoparametric hypersurfaces in the case $g=6$, only the first eigenvalue is computed to be the dimension $n$ by Tang-Yan \cite{TY13} and Tang-Xie-Yan \cite{TXY14} (and references therein, see also \cite{QT16} for more eigenvalue estimates), which solved Yau's conjecture on the first eigenvalue for isoparametric case: $\lambda_1(M^n)=n$ for any closed embedded minimal hypersurface $M$ in $\mathbb{S}^{n+1}$.

In view of the Chern Conjecture above, one can consider Yau's conjecture for minimal CSC hypersurfaces. Combining with Chang's result, we see that it is true in dimension $3$.

Solomon \cite{So92} showed that $2n$ is an eigenvalue for all $g=4$ minimal isoparametric hypersurfaces with certain quadratic forms (which vanish on one focal submanifold) being the eigenfunctions. Thus $2n$ could possibly be the second eigenvalue in the $g=4$ case. Interestingly, $3n$ is also an eigenvalue with certain cubic forms being the eigenfunctions (e.g., $\langle\nabla F(x)-4c_0|x|^2x, a\rangle$ for any $a\in\mathbb{R}^{n+2}$ where $F(x)$ is the Cartan-M\"{u}nzner polynomial and $c_0=(m_{-}-m_+)/(m_-+m_+)$ corresponds to the minimal isoparametric hypersurface $M_{c_0}:=F^{-1}(c_0)\cap \mathbb{S}^{n+1}$).

A related problem is Solomon's question on whether both focal submanifolds of $g=4$ isoparametric hypersurfaces are quadratic varieties. It has been shown in each isoparametric foliation one of the focal submanifolds, say $M_+$, is quadratic, which then provides eigenfunctions of the (possibly) second eigenvalue $2n$ for the minimal isoparametric hypersurface as mentioned above. However, from the classification of the nonnegative isoparametric polynomials $G_F^+(x):=|x|^4+F(x)$ $^\dag$\footnote{$^\dag$ nonnegative because of $|\nabla f|^2=g^2(1-f^2)$ by (\ref{CMeq}), where $f=F|_{\mathbb{S}^{n+1}}$ and $F(x)$ is any Cartan-M\"{u}nzner polynomial of $g=4$} which are sums of squares of quadratic forms, in many cases the other focal submanifold $M_-$, as the spherical zeroes of $G_F^+(x)$, are shown to be non-quadratic (see \cite{GT18}). In fact, it is shown that any quadratic form vanishing on these $M_-$ is identically zero. On the other hand, in the left cases there are indeed (non-zero) quadratic forms vanishing on $M_-$ though $G_F^+(x)$ is not a sum of squares of quadratic forms, thus leaving Solomon's question mysterious in these left cases.

\begin{problem}
Classify smooth minimal hypersurfaces asymptotic to hypercones over minimal isoparametric hypersurfaces of $\mathbb{S}^{n+1}$ in $\mathbb{R}^{n+2}$.
\end{problem}

For $g=1$, the hypercone is just a hyperplane and so is the asymptotic minimal hypersurface. This was firstly shown for area-minimizing hypersurfaces by Fleming \cite{Fl62} in his famous proof of the Bernstein Theorem, and for general minimal hypersurface one could use the monotonicity and the characterization of the density at infinity:
$\theta_{\infty}(\Sigma)=1$ if and only if $\Sigma$ is a hyperplane.

For $g=2$, the hypersones are just the Simons cones $C^{k,n-k}:=C(\mathbb{S}^k(\sqrt{\frac{k}{n}})\times \mathbb{S}^{n-k}(\sqrt{\frac{n-k}{n}}))$ (cones over the minimal Clifford tori). The asymptotic minimal hypersurfaces were classified for minimizing Simons cones ($n>6$, and $1<k<4$ if $n=6$) by Simon-Solomon \cite{SS86} and in general by Mazet \cite{Ma17}. There are exactly two (one when $n=2k$) asymptotic smooth minimal hypersurfaces up to similarity. Recall that in either components of $\mathbb{R}^{n+2}\setminus C$ there exists a homothety-invariant foliations by smooth minimizing hypersurfaces asymptotic to $C$, provided that $C$ is an area-minimizing cone (cf. \cite{SS86}, for $n=2k$ by Bombieri-De Giorgi-Giusti, $n\neq 2k$ by Simoes, and for general area-minimizing cones by Hardt-Simon \cite{HS85}).

For $g\geq3$ (i.e., $g=3,4,6$), the methods of Simon-Solomon \cite{SS86} and Mazet \cite{Ma17} failed to give such uniqueness. For minimizing ones (i.e., when $n+2\geq4g$ or $(g,m_-,m_+)\neq(4,1,6)$, cf. \cite{TZ19}) in these cases, White \cite{Wh89} obtained a third asymptotic minimizing hypersurface, but it might be non-smooth and admit some singularity as mentioned by himself. Therefore, it still leaves us the expectation of the uniqueness as in the $g\leq2$ cases above.

A related problem is Solomon-Yau's conjecture (cf. \cite{Yau90}): the lowest density among all non-hyperplane area-minimizing hypercones is attained by the Simons cones $C^{[\frac{n}{2}], [\frac{n+1}{2}]}$. Ilmanen-White \cite{IW15} showed that this is true in the class of topologically nontrivial cones in the asymptotic sense: the density is greater than $\sqrt{2}=\lim\limits_{m\rightarrow\infty}\theta_{\infty}(C^{m,m})$. The proof of the Willmore conjecture by Marques-Neves \cite{MN14} implies that the Solomon-Yau conjecture is true for $n=2$.

\begin{problem}
Scalar and Ricci curvature rigidity of isoparametric tubes in $\mathbb{S}^{n+1}$.
\end{problem}

Let $F(x)$ be a Cartan-M\"{u}nzner polynomial of degree $g$ on $\mathbb{R}^{n+2}$, $f:=F|_{\mathbb{S}^{n+1}}$ and $M_{c_0}:=f^{-1}(c_0)$ ($c_0:=(m_{-}-m_+)/(m_-+m_+)$) be the unique minimal isoparametric hypersurface in the isoparametric foliation $\mathfrak{F}:=\{M_t:=f^{-1}(t)\mid t\in [-1,1]\}$. Then $M_{c_0}$ divides the sphere into two connected components, say
\begin{equation}\label{DDBDSphere}
\mathbb{S}^{n+1}=\mathbb{S}^{n+1}_{F-}\bigcup_{M_{c_0}} \mathbb{S}^{n+1}_{F+}, \quad \mathbb{S}^{n+1}_{F-}:=\bigcup_{-1\leq t\leq c_0}M_t,\quad \mathbb{S}^{n+1}_{F+}:=\bigcup_{c_0\leq t\leq1}M_t,
 \end{equation}
 which are exactly tubes of radius $\frac{\pi}{g}-\theta_0$ around the focal submanifold $M_-=f^{-1}(-1)$ and of radius $\theta_0:=\frac{1}{g}\arccos(c_0)\in(0, \frac{\pi}{g})$ around $M_+=f^{-1}(1)$  respectively. Replacing $F$ by $-F$ will replace the sign $\pm$ above and thus we can henceforth consider without loss of generality the problem on one tube around $M_+$. For $c\in (-1, 1)$ and correspondingly $\theta:=\frac{1}{g}\arccos(c)\in (0, \frac{\pi}{g})$, we also denote the tube of radius $\theta$ around $M_+$ by
 \begin{equation}\label{isoptube}
 \mathcal{T}_{F+\theta}=\mathbb{S}^{n+1}_{F+c}=\bigcup_{c\leq t\leq1}M_t.
 \end{equation}

 For a Riemannian metric $\mathbf{g}$ on a manifold $N$ with boundary $\partial N$, we denote by $R_{\mathbf{g}}, \Ric_{\mathbf{g}}, H_{\mathbf{g}}, \II_{\mathbf{g}}$ the scalar curvature, the Ricci curvature on $N$, the mean curvature and the second fundamental form of $\partial N$ with respect to the inward unit normal vector field.
 Let $\bar{\mathbf{g}}$ denote the standard metric on the unit sphere.

 \begin{conjecture}\label{conj-scal}
 Let $\mathbf{g}$ be a Riemannian metric on the isoparametric tube $\mathbb{S}^{n+1}_{F+c}$ of $(\ref{isoptube})$ for $c\in (c_0,1)$ with the following properties:
 \begin{itemize}
 \item[$\bullet$] $R_{\mathbf{g}}\geq n(n+1)$ at each point in $\mathbb{S}^{n+1}_{F+c}$.
 \item[$\bullet$] The metrics $\mathbf{g}$ and $\bar{\mathbf{g}}$ induce the same metric on $M_c=\partial \mathbb{S}^{n+1}_{F+c}$.
 \item[$\bullet$] $H_{\mathbf{g}}\geq H_{\bar{\mathbf{g}}}$ at each point on $M_c$.
  \end{itemize}
 If $\mathbf{g}-\bar{\mathbf{g}}$ is sufficiently small in the $C^2$-norm, then $\varphi^*(\mathbf{g})=\bar{\mathbf{g}}$ for some diffeomorphism $\varphi: \mathbb{S}^{n+1}_{F+c}\rightarrow \mathbb{S}^{n+1}_{F+c}$ with $\varphi|_{M_c}=\id$.
 \end{conjecture}

 This is motivated by the scalar curvature rigidity derived from the positive mass theorem (cf. \cite{Br12}).

 For the $g=1$ case, the isoparametric tube is just a geodesic ball inside the hemisphere ($c_0=0$ now), and the conjecture holds for $c\geq\frac{2}{\sqrt{n+4}}$  by Brendle-Marques \cite{BM11}. As remarked in  \cite{BM11}, one cannot expect the same conclusion for the tube $\mathbb{S}^{n+1}_{F+}$ with $c=c_0$ (hemisphere in the $g=1$ case), because of the counterexamples to Min-Oo's conjecture constructed by Brendle-Marques-Neves \cite{BMN11}. However, the geodesic ball with $c\geq\frac{2}{\sqrt{n+4}}$ above can be enlarged slightly as shown in \cite{CMT13}. Therefore, it is convincing to have the Brendle-Marques' type scalar curvature rigidity on the isoparametric tube $\mathbb{S}^{n+1}_{F+c}$ with $c>c_0$ (i.e., of the ``biggest" radius $\theta<\theta_0$).

 Notice also that, because of the counterexamples to Min-Oo's conjecture in \cite{BMN11}, one cannot omit the assumption of the sufficiently smallness of $\mathbf{g}-\bar{\mathbf{g}}$ directly from the conjecture, unless one consider the rigidity on the isoparametric tube $\mathcal{T}_{F+\theta}$ of sufficiently small radius $\theta$. For this consideration, one can even assume $\II_{\mathbf{g}}\geq \II_{\bar{\mathbf{g}}}$ in place of $H_{\mathbf{g}}\geq H_{\bar{\mathbf{g}}}$ on $M_c$, in light of Min-Oo's conjecture and the following Ricci curvature rigidity.

  \begin{conjecture}\label{conj-Ric}
 Let $(N^{n+1}, \mathbf{g})$ be a compact Riemannian manifold with boundary $\partial N$ satisfying the following properties:
 \begin{itemize}
 \item[$\bullet$] $\Ric_{\mathbf{g}}\geq n\mathbf{g}$ in $N$.
 \item[$\bullet$] $(\partial N, \mathbf{g}|_{\partial N})$ is isometric to the isoparametric hypersurface $(M_c, \bar{\mathbf{g}}|_{M_c})$ for some $c\in [c_0,1)$.
 \item[$\bullet$] $\II_{\mathbf{g}}\geq \II_{\bar{\mathbf{g}}}$ at each point on $\partial N\cong M_c$.
  \end{itemize}
Then $(N^{n+1}, \mathbf{g})$ is isometric to the isoparametric tube $(\mathbb{S}^{n+1}_{F+c}, \bar{\mathbf{g}})$ in $\mathbb{S}^{n+1}$.
 \end{conjecture}

 The hemisphere case ($g=1$) was verified by Hang-Wang \cite{HW09} using Reilly's formula. This was generalized by Miao-Wang \cite{MW16} to domains of hemisphere with the second assumption replaced by an isometric immersion of $\partial N$ into unit sphere of possibly high codimension and the third assumption by $\II_{\mathbf{g}}\geq |\II_{\bar{\mathbf{g}}}|$, the latter of which forces the convexity of the boundary $\partial N$ and thus rules out isoparametric tubes with $g\geq2$ (which are always non-convex).

 Using a similar method as Hang-Wang, Lai \cite{La19} recently verified the $g=2$ case of Conjecture \ref{conj-Ric} (the isoparametric tubes bounded by the Clifford tori), with a slightly weaker assumption on the second fundamental form.

 \begin{problem}
Serrin's overdetermined problem on isoparametric tubes in $\mathbb{S}^{n+1}$.
\end{problem}

  \begin{conjecture}\label{conj-Ser}
  Let $\Omega\subseteq\mathbb{S}^{n+1}$ be a domain with smooth connected boundary $\partial \Omega$.
  Suppose there exists a positive solution $u\in C^2(\overline{\Omega})$ of Serrin's overdetermined problem:
  \begin{equation} \label{equ-ser}
  \left\{\begin{array}{ll}
  \Delta u = a(u), & in~ \Omega\\
  u=0,~ \frac{\partial u}{\partial \nu}=\const & on~\partial\Omega,
  \end{array}
  \right.
  \end{equation}
  where $a\in C^1(\mathbb{R})$ and $\frac{\partial u}{\partial \nu}$ denotes the normal derivative along the inward normal field on $\partial\Omega$.
  Then
  \begin{itemize}
  \item[$(i)$] $\Omega=\mathbb{S}^{n+1}_{F+c}$ is an isoparametric tube for some $F(x)$ and $c\in (-1, 1)$;
  \item[$(ii)$] $u$ is constant on $M_t$ for each $t\in [c,1]$, i.e., $u$ is a function of $f=F|_{\mathbb{S}^{n+1}}$.
  \end{itemize}
  \end{conjecture}

  Serrin's classical result \cite{Se71} asserts that if $\Omega\subset \mathbb{R}^{n+1}$ is a bounded and smooth domain for which there is a positive solution to the equation (\ref{equ-ser}) then $\Omega$ is a ball and $u$ is radially symmetric (exactly the only bounded isoparametric tube in the Euclidean space). For unbounded domains of $\mathbb{R}^{n+1}$, there are many studies and developments on Serrin's overdetermined problem, e.g., the Berestycki-Caffarelli-Nirenberg conjecture which has the same assertion as Conjecture \ref{conj-Ser} that $\Omega$ is an isoparametric tube, namely a half-space or a cylinder in this case. However, there are several counterexamples and also affirmative answers under certain additional assumptions. We refer to the survey \cite{CW18} and references therein for more details in this direction.

 As Alexandrov's theorem, Serrin's classical result has also been generalized successfully to domains in the hemisphere and the hyperbolic space (cf. \cite{KP98}). Also similarly, in the whole sphere we have more domains like isoparametric tubes which support solutions to the Serrin's overdetermined equation (\ref{equ-ser}). However, if without the connectedness assumption of $\partial\Omega$, then a non-isoparametric domain, in fact a neighborhood of the great hypersphere $\mathbb{S}^n\subset\mathbb{S}^{n+1}$, has been constructed to admit solutions to (\ref{equ-ser}) by Fall-Minlend-Weth \cite{FMW18}. We refer to \cite{Sa18, Sa16} for more related discussions, where isoparametric tubes in unit spheres are characterized by the so-called constant flow property of heat flow.

 Comparing with the hemisphere case above, one can assume further in Conjecture \ref{conj-Ser} that $\Omega$ is a domain in a given isoparametric tube $\mathbb{S}^{n+1}_{F+}$. However, unlike the hemisphere case, we warn that the final isoparametric tube in the assertion could be of different type as the given one, since isoparametric tubes with $g\geq2$ are ``big" enough to contain small geodesic balls or other smaller isoparametric tubes.

 Another motivation comes from a result of Wang \cite{Wang87} which states that the level hypersurfaces of a transnormal function on $\mathbb{R}^{n+1}$ or $\mathbb{S}^{n+1}$ are isoparametric (incorrect on hyperbolic space $\mathbb{H}^{n+1}$). This means that a solution to the first equation of (\ref{isop-def}) is sufficient to generate isoparametric hypersurfaces in $\mathbb{R}^{n+1}$ or $\mathbb{S}^{n+1}$. Now Conjecture \ref{conj-Ser} attempts to deduce the same conclusion from the second equation of (\ref{isop-def}) necessarily with certain boundary assumptions.


\section*{3. Problems related to isoparametric theory generally}

\begin{problem}
Classify isoparametric foliations on homotopy spheres.
\end{problem}

From the celebrated solutions of the (generalized) Poincar\'{e} conjectures by Smale, Freedman and Perelman, we know that homotopy spheres are homeomorphic to the standard spheres.
Due to the surprising discovery of exotic spheres by Milnor, we know that there are many homotopy spheres which are not diffeomorphic to the standard spheres.
On general exotic spheres we do not know good and canonical metrics (cf. \cite{JW08}). Therefore, to study isoparametric foliations on homotopy spheres we prefer to use the viewpoint of differential topology based on the Qian-Tang fundamental construction (see Section 1 or \cite{QT15}), which states that if the manifold $N=D(E_-)\bigcup_{\varphi}D(E_+)$ is glued from two disk bundles $D(E_{\pm})$ of rank greater than $1$ along their boundaries $\varphi\in \Diff(\partial E_-\rightarrow \partial E_+)$, then there exists a metric $\mathbf{g}$ such that the canonical singular foliation $\mathfrak{F}$ is an isoparametric foliation under this metric. We call such a manifold $N$ to have a \textit{DDBD} structure $(D(E_{\pm}), \varphi)$, which is then equivalent to saying that $N$ has an isoparamtric foliation when $N$ is a closed simply connected manifold (cf. \cite{QG16}).

As the famous result of Smale which says that every homotopy (exotic) sphere $\Sigma^n=D^n_-\bigcup_\phi D^n_+$ $(n\geq5)$ is a twisted sphere (glued from two disks) for some $\phi\in \Diff(\mathbb{S}^{n-1})$, we \cite{Ge16} showed the analogous result that $\Sigma^n=D(E_-)\bigcup_{\psi}D(E_+)$ is always a DDBD-twisted sphere if $\mathbb{S}^n=D(E_-)\bigcup_{\varphi}D(E_+)$. This immediately leads to the corollary that (equivalence classes of) isoparametric foliations on $\Sigma^n$ are in one-to-one correspondence with that on $\mathbb{S}^n$ $(n\geq5)$. Henceforth we shall need only
consider the classification of DDBD structures  $(D(E_{\pm}), \varphi)$ on $\mathbb{S}^n$ $(n\geq5)$. In dimension four, we \cite{GT13} have shown that an exotic $4$-sphere (if exist) does not admit any isoparametric foliations. On $\mathbb{S}^n$ with $n\leq4$, it is known that there are only the standard ones (cf. \cite{GT13, Ra12}).

In \cite{Ge16} it has been shown that isoparametric foliations $(\mathbb{S}^n, \mathfrak{F})$ with two points as the focal submanifolds are unique up to foliated diffeomorphisms for $n\neq 5$, and are of the same number as that of components of $\Diff^+(\mathbb{S}^4)$ for $n=5$ which is still unknown so far.

Notice that there are more isoparametric foliations (DDBD structures) found on $\mathbb{S}^n$ than the classical isoparametric foliations under the round metric (cf. \cite{Ge16}), each of which have homeomorphic isoparametric hypersurfaces and focal submanifolds. An interesting question of Fang \cite{Fa17} asks that whether this homeomorphic phenomenon holds in general.

\begin{problem}
Existence of isoparametric functions on homotopy spheres.
\end{problem}

 Under a given Riemannian metric on a homotopy sphere $\Sigma^n$, it is still a widely unknown and very hard problem that whether there exist a properly isoparametric (or just properly transnormal) function, namely a nonconstant solution to equation (\ref{isop-def}) with connected level sets (see the notion of properly isoparametric in \cite{GT13}). It worths to remark that a positive answer to this question in dimension $n=4$ will solve affirmatively the $4$ dimensional smooth Poincar\'{e} conjecture: there is no exotic $4$-sphere, by the result of \cite{GT13} mentioned in last problem. In fact it is sufficient to show the existence of a properly transnormal function or even a singular Riemannian foliation on a $4$-sphere $\Sigma^4$ under some given Riemannian metric (cf. \cite{GT13, GR15}).  Certain geometric flows are expected to be helpful to this question.

\begin{problem}
Compact version of Cheeger-Gromoll's soul theorem.
\end{problem}

\begin{conjecture}\label{Groveconj}
$(\GroveConj)$ A closed simply connected manifold with nonnegative sectional curvature has a DDBD structure, i.e., an isoparametric foliation under some metric (cf. \cite{GR15}).
\end{conjecture}

Cheeger-Gromoll's soul theorem states that a complete open manifold with $K_{sec}\geq0$ is diffeomorphic to the normal bundle of a compact submanifold (called \textit{soul}) of $K_{sec}\geq0$.
Therefore, the Grove conjecture above is rather natural by replacing one vector bundle for open manifold with two vector bundles for compact manifold.
An immediate but important application will be a differential classification of all closed simply connected $4$-manifolds with nonnegative sectional curvature, namely the standard $\mathbb{S}^4$, $\mathbb{CP}^2$, $\mathbb{S}^2\times \mathbb{S}^2$ and $\mathbb{CP}^2\#\pm \mathbb{CP}^2$ (cf. \cite{GR15}).

In other words, this question asks whether nonnegative curvature implies isoparametric foliation. In the following we ask the converse under leafwise nonnegativity.

\begin{problem}
Does the leafwise nonnegativity imply the whole nonnegativity of curvatures on an isoparametric foliation $(N,\mathfrak{F})$?
\end{problem}
\begin{conjecture}\label{conj-leafwise}
Given an isoparametric foliation $(N^{n}, \mathfrak{F})$ with nonnegative or positive leafwise sectional curvature $K_{sec}^F\geq0$ or $K_{sec}^F>0$ (resp. nonnegative/positive leafwise Ricci curvature $\Ric^F$ or leafwise scalar curvature $R^F$), then the whole manifold $N$ admits a metric of nonnegative or positive curvatures correspondingly.
\end{conjecture}

It is satisfactory to consider the conjecture only on closed simply connected manifolds. Then there are exactly two singular leaves, the focal submanifolds $M_\pm$. Without loss of generality we can assume the leaf space to be the interval $[-1, 1]$ and then $$\left(N':=N\setminus (M_-\cup M_+), \quad \mathfrak{F}':=\mathfrak{F}|_{N'}=\{M_t \mid M_t\cong M_0, t\in (-1,1)\}\right)$$ is a regular Riemannian foliation of codimension one whose integral subbundle is denoted by $F'\subset TN'$ .

A one-parameter family of smooth metrics $\mathbf{g}_t$, $t\in [-1,1]$, on $M_t$ is called a family of \textit{smooth leafwise metrics} on a codimension one singular foliation $(N,\mathfrak{F})$ if $\mathbf{g}_t$ is the restriction to $M_t$ of a foliated metric $\mathbf{g}^F$ on $(N,\mathfrak{F})$. By a \textit{foliated metric} $\mathbf{g}^F$ on a singular foliation $(N,\mathfrak{F})$ we mean it is a smooth symmetric positive definite $(2,0)$-tensor generated by its valuation on the vector fields $X_i's$ that span the tangent bundles of the leaves. It turns out that the restriction of smooth leafwise metrics $\mathbf{g}_t$ to $t\in (-1,1)$ equals some Euclidean metric $\mathbf{g}^{F'}$ on the integral subbundle $F'$ over the regular foliation $(N', \mathfrak{F}')$. The restriction of the bundle-like Riemannian metric $\mathbf{g}^N$ of the singular Riemannian foliation is naturally a foliated metric and thus provides the canonical family of smooth leafwise metrics $\mathbf{g}^N_t$.

Given a family of smooth leafwise metric $\mathbf{g}_t$ $(t\in [-1,1])$,  we can then define the leafwise sectional, Ricci and scalar curvatures $K_{t}^F, \Ric^F_t, R^F_t$ to be the curvatures defined by the Riemannian metric $\mathbf{g}_t$ on $M_t$; and call that $(N^{n+1}, \mathfrak{F})$ has leafwise nonnegative or positive curvature if the leafwise curvature is so for all  $t\in [-1,1]$, denoted as in the conjecture above. We remark that these definitions are also suitable to general (singular) foliations with one-parameter $t$ replaced by the leaf space $t\in N/\mathfrak{F}$.

A beautiful result of Zhang \cite{Zhang17} states that if a regular foliation $(N^n, \mathfrak{F})$ $(n\geq5)$ has positive leafwise scalar curvature $R^F>0$ then $N$ admits a Riemannian metric of positive scalar curvature. Hence the conjecture above generalizes this to the case of singular foliations but under some Riemannian conditions (SRF or isoparametric) which one would also boldly omit in comparison with Zhang's result. On the other hand, the condition of isoparametric foliations alone does not imply any positivity of the curvature, since the remarkable result of Hitchin \cite{Hi74} shows that there are exotic spheres $\Sigma^n$ $(n\geq9)$ admitting no metrics of positive scalar curvature, while $\Sigma^n$ indeed supports isoparametric foliations as mentioned before.

Also due to Hitchin's result, the following interesting problem of Tang seems contradictory to Conjecture \ref{conj-leafwise}:
\begin{conjecture}
$(\TangConj)$ Every isoparametric hypersurface and focal submanifold of $\mathbb{S}^{n}$ admit a metric of nonnegative sectional curvature.
\end{conjecture}
This conjecture can be extended to homotopy spheres (and thus Hitchin's exotic spheres $\Sigma^n$ with no metrics of positive scalar curvature) by the one-to-one correspondence of \cite{Ge16} if we ignore the round metric as before. For those homogeneous isoparametric foliations, the Tang conjecture holds true naturally as all leafs would admit a normal homogeneous metric of nonnegative sectional curvature.

Here we explain why we think the Tang conjecture does not contradict Conjecture \ref{conj-leafwise}. Briefly speaking, it would happen that the non-negatively curved metrics on the isoparametric hypersurfaces and the focal submanifolds would not constitute a family of smooth leafwise metrics $\mathbf{g}_t$ on the whole foliation $(\Sigma^n, \mathfrak{F})$.

\end{document}